\documentclass[12pt,reqno]{amsart}

\newtheorem{theorem}{Theorem}[section]
\newtheorem{lemma}[theorem]{Lemma}

\theoremstyle{definition}

\theoremstyle{remark}

\numberwithin{equation}{section}

\usepackage{times}
\usepackage{enumerate}
\usepackage{graphicx,adjustbox}
\usepackage{hyperref}
\usepackage{amsmath,amssymb}

\newcommand{\lcm}{{\rm \mbox{lcm}}}

\usepackage{amscd}
\usepackage{graphicx}
\usepackage[all]{xy}

\title[Affine monomial curves]
{Affine monomial curves}
\author{
Indranath Sengupta
}
\date{}

\address{\small \rm  Discipline of Mathematics, IIT Gandhinagar, Palaj, Gandhinagar, 
Gujarat 382355, INDIA.}
\email{indranathsg@iitgn.ac.in}
\thanks{This author thanks SERB for their support through 
the MATRICS grant MTR/2018/000420.}

\date{}



\allowdisplaybreaks

\begin{document}

\begin{abstract}
This article is an expository survey on affine monomial curves, where 
we discuss some research problems from the perspective of computation.
\end{abstract}

\maketitle

This article is an attempt to introduce and discuss some ongoing 
research in the field of affine monomial curves and numerical 
semigroup rings. There is a large body of literature in this 
area and we would not make any attempt to cite most of them. 
We would rather cite some, which have a strong list of 
references, \textit{e.g.}, \cite{rgs} and \cite{stamate}. The 
topics which we have chosen to discuss for this article are limited 
and should not be misinterpreted as an exhaustive list. We have 
chosen the topics from the perspective of computation, where computer 
algebra software like Singular \cite{singular} and GAP \cite{GAP4} can 
be used effectively.  This article could be useful for some young 
researchers, who are interested in the computational aspects 
of commutative algebra and algebraic geometry. The texts 
\cite{matsumura}, \cite{bh}, \cite{peeva} and \cite{iyengar} provide 
most of the background material in commutative algebra and 
\cite{singularbook}, \cite{dl} are strongly recommended for learning 
the subject using computer algebra as a companion.

\section{Numerical Semigroups and Affine Monomial Curves}

Let $\mathbb{N}$ denote the set of nonnegative integers. 
A \textit{numerical semigroup} $\Gamma$ is a subset of 
$\mathbb{N}$ containing $0$, closed under addition and 
generates $\mathbb{Z}$ as a group. It follows that 
(see \cite{rgs}) the set $\mathbb{N}\setminus \Gamma$ is 
finite and that the semigroup $\Gamma$ has a unique minimal 
system of generators $n_{0} < n_{1} < \cdots < n_{p}$, that is, 
$$\Gamma=\Gamma(n_0,\ldots, n_p) = 
\lbrace\sum_{j=0}^{p}z_{j}n_{j}\mid z_{j}\quad \text{nonnegative \, integers}\rbrace.$$ 
The term minimal stands for the following: if $n_i=\sum_{j=0}^{p}z_{j}n_{j}$ 
for some non-negative integers $z_{j}$, then $z_{j}=0$ for 
all $j\neq i$ and $z_{i}=1$. 
The integers $n_{0}$ and $p + 1$ are known 
as the \textit{multiplicity} $m(\Gamma)$ and the 
\textit{embedding dimension} $e(\Gamma)$. The greatest 
integer not in $\Gamma$ is called the 
\textit{Frobenius number} of $\Gamma$, denoted by 
$F(\Gamma)$. The integer $c(\Gamma) = F(\Gamma) + 1$ is called 
the \textit{condcutor} of $\Gamma$. The \textit{Ap\'{e}ry set} 
of $\Gamma$ with respect to a non-zero $a\in \Gamma$ is defined 
to be the set $\rm{Ap}(\Gamma,a)=\{s\in \Gamma\mid s-a\notin \Gamma\}$. 
The numerical semigroup $\Gamma$ is \textit{symmetric} if 
$F(\Gamma)$ is odd and $x\in \mathbb{Z}\setminus \Gamma$ 
implies $F(\Gamma)-x\in \Gamma$. For example, if $\Gamma$ is the 
numerical semigroup generated minimally by the integers 
$3, 5, 7$, then $m(\Gamma)=3$, $e(\Gamma) = 3$, $F(\Gamma) = 4$, 
$c(\Gamma)=5$. It is easy to verify that $\Gamma$ is not a symmetric 
semigroup. Computing $F(\Gamma)$ is a hard problem in general, 
known as the \textit{Frobenius Problem} (see \cite{alfonsin}) or the 
\textit{coin exchange problem} (see \cite{ian}) . 
\medskip

Let $p\geq 1$ and $n_{0}, \ldots, n_{p}$ be positive integers with 
$\gcd (n_{0},\ldots,\, n_{p})=1$. Let us assume that the numbers 
$n_{0}, \ldots, n_{p}$ generate the numerical semigroup 
$\Gamma(n_0,\ldots, n_p)$ minimally. Let $k$ denote a field and 
$\eta:k[x_0,\,\ldots,\, x_p]\rightarrow k[t]$ be the mapping defined by 
$\eta(x_i)=t^{n_i},\,0\leq i\leq p$. Let $\frak{p}(n_0,\ldots, n_p) = 
\ker (\eta)$. The map $\eta$ is a monomial parametrization and defines 
a curve $\mathcal{C} = \{(t^{n_{0}}, \ldots, t^{n_{p}}) \mid t\in k\}$, 
known as a \textit{monomial curve} in the affine space $\mathbb{A}^{p+1}$. 
The ideal $\frak{p}(n_0,\ldots, n_p)$ is called the defining ideal of the 
curve $\mathcal{C}$, which has been called by the name \textit{monomial primes} 
by some authors; \cite{cutkosky}, \cite{hema}. The affine $k$-algebra 
$A = k[x_0,\,\ldots,\, x_p]/\frak{p}(n_0,\ldots, n_p)=k[t^{n_{0}}, \ldots ,t^{n_{p}}]$ 
is called the \textit{coordinate ring} of the affine monomial curve $\mathcal{C}$. 
Let $R$ denote the polynomial ring $k[x_0,\,\ldots,\, x_p]$. The ideal 
$\frak{p}(n_0,\ldots, n_p)$ is a prime ideal with dimension one, that is, 
the Krull dimension of the affine $k$-algebra $A$ is one, therefore, $A$ is 
Cohen-Macaulay and $\frak{p}(n_0,\ldots, n_p)$ is a perfect ideal.

\section{Defining Equations, Gr\"{o}bner bases and Syzygies}
The ideal $\frak{p}(n_0,\ldots, n_p)$ is a graded ideal with respect 
to the weighted gradation given by ${\rm wt}(x_{i})=n_{i}$, ${\rm wt}(t)=1$ 
on the polynomial rings $k[x_0,\,\ldots,\, x_p]$ and $k[t]$. Therefore, 
by the graded version of Nakayama's lemma, the minimal number of 
generators of $\frak{p}(n_0,\ldots, n_p)$ is a well-defined notion 
and this number is denoted by $\beta_{1}(\frak{p}(n_0,\ldots, n_p))$, 
known as the \textit{first Betti number} of $\frak{p}(n_0,\ldots, n_p)$. 
It is both interesting and hard to compute a minimal generating set and 
$\beta_{1}(\frak{p}(n_0,\ldots, n_p))$ for the defining ideal 
$\frak{p}(n_0,\ldots, n_p)$. It is not difficult to show that the ideal 
$\frak{p}(n_0,\ldots, n_p)$ is generated by binomials; see Lemma 4.1 in 
\cite{sturmfels} for a more general statement on toric ideals). Choosing a finite 
number of them minimally usually requires a good understanding of the 
Ap\'{e}ry set; see \cite{patsingh}, \cite{patil2}, \cite{malooseng}, \cite{mss1}. The 
computer algebra package GAP \cite{GAP4} is quite a useful companion for 
these calculations. It is 
usually harder to compute a generating set for 
$\frak{p}(n_0,\ldots, n_p)$, which forms a Gr\"{o}bner basis. Let 
us recall the basics of Gr\"{o}bner basis first so that we can understand 
its usefulness. We refer to \cite{clo} for a good introduction to the 
subject.
\medskip

Let $k[x_{1},\ldots,x_{n}]$ be the polynomial ring over $k$, 
with a monomial order $>$. For $0 \neq f\in k[x_1,\ldots,x_n]$, 
let ${\rm Lm}_{>}(f)$, ${\rm Lt}_{>}(f)$, ${\rm Lc}_{>}(f)$ 
denote the {\it leading monomial}, {\it leading term} and 
{\it leading constant} of $f$ respectively. These are often 
written as ${\rm Lt}(f)$, ${\rm Lc}(f)$, ${\rm Lm}(f)$, omitting 
$``>"$ from the expression. For $0\neq f, g\in k[x_{1},\ldots,x_{n}]$, 
the {\it S - polynomial} $S(f, g)$, is the polynomial 
$$S(f, g) := \frac{\lcm({\rm Lm}(f), {\rm Lm}(g))}{{\rm Lt}(f)}\cdot f - 
\frac{\lcm({\rm Lm}(f), {\rm Lm}(g))}{{\rm Lt}(g)}\cdot g\,.$$ 
Given $G = \{g_{1},\ldots, g_{t}\}\subseteq k[x_1,\ldots,x_n]$ and 
$f\in k[x_{1},\ldots,x_{n}]$, we say that 
$f$ {\it reduces to zero modulo} $G$, 
denoted by $f\to_{G} 0$, if $f$ can be written as 
$f = \sum_{i=1}^{t} a_{i}g_{i}\,$, 
such that $\,{\rm Lm}(f) \geq {\rm Lm}(a_{i}g_{i})$, whenever 
$\,a_{i}g_{i}\neq 0$. Buchberger's Criterion says that, for an ideal $I$ 
in $k[x_{1},\ldots,x_{n}]$ generated by the set $G = \{g_{1},\ldots, g_{t}\}$, 
$G$ is a Gr\"obner basis for $I$ iff $S(g_{i}, g_{j}) \to_{G} 0$, for every $i\neq j$. 
The situation is simpler when $\gcd(\,{\rm Lm}(f), \,{\rm Lm}(g)\,) = 1$, and one 
chooses monomial orders in order to maximize the occurrence of such pairs in the list 
of polynomials $G = \{g_{1},\ldots, g_{t}\}$. The following Lemma is often useful in 
computing a Gr\"{o}bner basis.
\medskip

\begin{lemma} 
Let $G = \{g_{1},\ldots, g_{t}\}\subseteq k[x_{1},\ldots,x_{n}]$ and let 
$f, g\in G$ be non-zero with ${\rm Lc}(f) = {\rm Lc}(g) = 1$, and 
$\gcd(\,{\rm Lm}(f), \,{\rm Lm}(g)\,) = 1$. Then, 
\begin{enumerate}
\item $S(f, g) \,=\,  {\rm Lm}(g).\,f \,-\, {\rm Lm}(f).\,g$\,.

\item $S(f, g) \,=\,  -(g - {\rm Lm}(g)).\,f \,+\, (f - {\rm Lm}(f)).\,g\, 
\longrightarrow_{G} 0$. 

\end{enumerate}
\end{lemma}
\medskip

The computation with the $S$-polynomials can be used together 
with the Schreyer's theorem \ref{schreyer} for writing down the 
first syzygy for an ideal. We refer to the book by Peeva \cite{peeva} 
for syzygies and minimal free resolutions. Let us recall Schreyer's 
theorem from Chapter 5, Theorem 3.2 in \cite{cox}.

\begin{theorem}[Schreyer]\label{schreyer}{\it 
Let $K$ be a field, $k[x_{1},\ldots,x_{n}]$ 
be the polynomial ring and $I$ be an ideal in $k[x_{1},\ldots,x_{n}]$. 
Let $G := \{g_{1}, \ldots, g_{t}\}$ be an ordered set of generators 
for $I$, which is a Gr\"obner basis, with respect to some fixed monomial 
order on $k[x_{1},\ldots,x_{n}]$. Let 
${\rm Syz}(g_{1}, \ldots, g_{t}) := \{(a_{1}, \ldots, a_{t})\in R^{t} \mid 
\sum_{i=1}^{t}a_{i}g_{i} = 0\}$. 
Suppose that, for $\,i\neq j$, 
$$ S(g_{i}, g_{j})  =  a_{j}g_{i} - a_{i}g_{j} 
= \sum_{k=1}^{t}h_{k}g_{k}\, \longrightarrow_{G} 0\,.$$
Then, the 
$\,t$-tuples 
$$\left( \begin{array}{lclclcl}
h_{1}\,, & \cdots \,, & h_{i} - a_{j}\,,  & \cdots \,, 
& h_{j} + a_{i}\,,  & \cdots \,, & h_{t}
\end{array}\right)\,$$ 
generate $\,{\rm Syz}(g_{1}, \ldots, g_{t})$.
}
\end{theorem}
\medskip

The ideal $\mathfrak{p}(n_{0},\ldots , n_{p})$ is a graded ideal, the 
$R$-module $A$ is a graded module with projective dimension $p$. Therefore, 
the $R$-module $A$ has a minimal graded free resolution 
$$0 \longrightarrow R^{\beta_{p}} \stackrel{\vartheta_{p-1}}
\longrightarrow R^{\beta_{p-2}}\longrightarrow \cdots 
\longrightarrow R^{\beta_{2}}\stackrel{\vartheta_{2}}
\longrightarrow R^{\beta_{1}} 
\stackrel{\vartheta_{1}}\longrightarrow R 
\stackrel{\vartheta_{0}}\longrightarrow A \longrightarrow 0\,.$$ 
The Betti numbers $\beta_{i}$ will be denoted by 
$\beta_{i}(\mathfrak{p}(n_{0},\ldots , n_{p}))$, often called the 
\textit{Betti numbers} of $\mathfrak{p}(n_{0},\ldots , n_{p})$. These 
are important numerical invariants for the curve $\mathcal{C}$, 
particularly, for understanding an embedding of $\mathcal{C}$ in 
$\mathbb{A}_{k}^{p+1}$.
\medskip

One can see computation of Gr\"{o}bner basis and application of 
Schreyer's theorem in the articles \cite{seng1}, \cite{seng2} 
and \cite{rst}. This is often the only effective way to compute a 
free resolution for the ideal $\mathfrak{p}(n_{0},\ldots, n_{p})$, with 
some amount of help from a computer algebra software like 
Singular \cite{singular} in guessing a candidate for the Gr\"{o}bner 
basis of the ideal. Sometimes, one becomes lucky when the 
ideal $\mathfrak{p}(n_{0},\ldots, n_{p})$ has some special structure. 
For example, in the case when $n_{0}, \ldots , n_{p}$ is a minimal 
arithmetic sequence, the defining ideal has a special structure, 
which allows one to adopt some homological techniques like mapping 
cone to write a free resolution; see \cite{gss1}, \cite{gss2}. However, 
Schreyer's technique or mapping cone or other techniques rarely create a 
free resolution that is minimal. The next difficulty is to spot the 
redundancy in the free resolution obtained from one of the techniques 
and extract a minimal free resolution which sits as a direct summand of 
the free resolution; see \cite{gss2}, \cite{rst}. Some questions 
emerge from the above discussions and let us list them down.
\medskip

\begin{enumerate}
\item[Q1.] Compute a minimal generating set for the defining ideal 
$\mathfrak{p}(n_{0},\ldots, n_{p})$ and write a formula for the 
first Betti number $\beta_{1}(\frak{p}(n_0,\ldots, n_p))$.
\medskip

\item[Q2.] Compute a Gr\"{o}bner basis for the defining ideal 
$\mathfrak{p}(n_{0},\ldots, n_{p})$ and use Schreyer's theorem 
to compute the first syzygy module. Compute a minimal subset 
of the first syzygy module and write a formula for the cardinality 
of that set, known as the second Betti number, denoted by 
$\beta_{2}(\frak{p}(n_0,\ldots, n_p))$.
\medskip

\item[Q3.] Write a minimal free resolution for the defining ideal 
$\mathfrak{p}(n_{0},\ldots, n_{p})$ and compute all the Betti 
numbers $\beta_{i}(\frak{p}(n_0,\ldots, n_p))$.
\end{enumerate}

\section{Betti numbers and their unboundedness}
The $i$-th Betti number $\beta_{i}(\frak{p}(n_0,\ldots, n_p))$ 
of the ideal $\frak{p}(n_0,\ldots, n_p)$ is the rank of the $i$-th 
free module in a minimal free resolution of $\frak{p}(n_0,\ldots, n_p)$. 
For a given $p\geq 2$, let 
$\beta_{i}(r) = { \sup}(\beta_{i}(\frak{p}(n_1,\ldots, n_r)))$, 
where $\sup$ is taken over all the sequences of positive integers 
$n_0,\ldots, n_p$. Herzog \cite{herzog} proved that $\beta_{1}(3)=3$. 
It is not difficult to show that $\beta_{2}(3)$ is a finite integer 
as well. Bresinsky in \cite{bre}, \cite{bre1}, \cite{bre2}, \cite{bre3}, 
\cite{brehoa}, extensively studied relations among the 
generators $n_0,\ldots, n_p$ of the numerical semigroup defined by 
these integers. It was proved in  \cite{bre1} and \cite{bre2} respectively 
that, for $r=4$ and for certain cases 
in $r=5$, the symmetry condition on the semigroup generated by 
$n_0,\ldots, n_p$ imposes an upper bound on 
$\beta_{1}(\frak{p}(n_0,\ldots, n_p))$. The following question is 
open in general:
\medskip

\begin{enumerate}
\item[Q4.] Given $p\geq 2$, does the symmetry condition on the numerical 
semigroup minimally generated by $n_0,\ldots, n_p$ impose 
an upper bound on $\beta_{1}(\frak{p}(n_0,\ldots, n_p))$?
\end{enumerate}  
\medskip

Bresinsky \cite{bre} 
constructed a class of monomial curves in $\mathbb{A}^{4}$ to prove that 
$\beta_{1}(4)=\infty$. He used this observation to prove that 
$\beta_{1}(r)=\infty$, for every $r\geq 4$. Let us recall Bresinsky's 
example of monomial curves in $\mathbb{A}^{4}$, as 
defined in \cite{bre}. Let $q_2\geq 4$ be even, 
$q_{1} = q_{2}+1,\, d_1=q_{2}-1$. Set 
$n_{0}=q_{1}q_{2},\, n_{1}=q_{1}d_{1},\, n_{2}=q_{1}q_{2}+d_{1},\, 
n_{3}=q_{2}d_{1}$. Clearly $\gcd (n_1,\, n_2,\, n_3,\, n_4)=1$. Let us use the shorthand $\mathbf{\underline{n}}$ to denote 
Bresinsky's sequence of integers defined above. Bresinsky \cite{bre} proved 
that the set $A= A_{1}\cup A_{2}\cup \{g_1, g_2\}$ generates the ideal 
$\mathfrak{p}(n_1,\ldots, n_4)$, where $A_{1}=\{f_{\mu}| f_{\mu}=x_{1}^{\mu-1} x_{3}^{{q_2}-
\mu}-x_{2}^{{q_2}-\mu}x_{4}^{\mu+1},\quad 1\leq\mu\leq q_{2}\}$, 
$A_{2}=\{f| f=x_{1}^{\nu_{1}} x_{4}^{\nu_{4}}-x_{2}^{\mu_{2}}x_{3}^{\mu_{3}},\, \nu_{1},\ \mu_{3}< d_{1}\}$ and $g_{1}=x_{1}^{d_{1}}-{x_{2}}^{q_{2}}$, $g_{2}=x_{3}x_{4}-x_{2}x_{1}$. We have proved in 
\cite{mss1} that for Bresinsky's examples $\beta_{i}(4) = \infty$, 
for every $1\leq i\leq 3$. We have also described all the syzygies 
explicitly. We have in fact proved that the set $S=A_{1}\cup A_{2}^{'}\cup \{g_1, g_2\}$, where $A_{2}'=\{h_{m}\mid x_{1}^{m}x_{4}^{(q_{1}-m)}-x_{2}^{(q_{2}-m)}x_{3}^{m}, 1\leq m\leq q_{2}-2\}$ is a minimal generating set and also 
Gr\"{o}bner basis for the ideal $\mathfrak{p}(\mathbf{\underline{n}})$ 
with respect to the lexicographic monomial order induced by $x_{3}>x_{2}>x_{1}>x_{4}$ on $k[x_{1}, \ldots , x_{4}]$. It turns out that 
$\beta_{1}(\mathfrak{p}(\mathbf{\underline{n}})) = \mid S \mid  = 2q_{2}$. 
Then it has been proved using Schreyer's theorem that 
$\beta_{2}(\mathfrak{p}(\mathbf{\underline{n}}))= 4(q_{2}-1)$ and 
also $\beta_{3}(\mathfrak{p}(\mathbf{\underline{n}}))= 2q_{2}-3$. 
A similar study has been carried out by J. Herzog and D.I. Stamate 
in \cite{herstamate} and \cite{stamate}. However, the objective and 
approach in our study are quite different. 
\medskip

The integers $n_{0}, n_{1}, n_{2}, n_{3}$ defining Bresinsky's 
semigroups have the property that $n_{0} + n_{3} = n_{1} + n_{2}$, 
up to a renaming if necessary. We generalize this construction in 
arbitrary embedding dimension by considering the family of semigroups 
$\Gamma(M)$, minimally generated by the set 
$M=\{a, a+d, a+2d,\ldots,a+(e-3)d, b,b+d\}$, $b > a+(e-3)d$, 
$\gcd(a,d)=1$ and $d\nmid(b-a)$. Note that the sequence of 
integers is a \textit{concatenation of two arithmetic sequences} and 
the condition $d\nmid(b-a)$ ensures that the sequence is not 
a part of an arithmetic sequence. See \cite{mss2} for 
discussion on this class of semigroups. In \cite{mss2} the case $a=e+1$ 
has been studies in detail and it has been proved that $\mu(\mathfrak{p})$ 
is bounded. In \cite{mss2} it has been shown that there are symmetric 
semigroups of this form, generalizing Roasales' construction in \cite{ros}. 
Our observations in \cite{mss2} give rise to the following 
question:
\medskip

\begin{enumerate}
\item[Q5.] Let $p\geq 3$ be an integer, $b > a+(p-2)d$, $\gcd(a,d)=1$ 
and $d\nmid(b-a)$. Let $M=\{a, a+d, a+2d,\ldots,a+(p-2)d, b,b+d\}$ be 
the set of $p+1$ positive integers. Let $\Gamma(M)=\langle M\rangle$ 
be the numerical semigroup generated by the set $M$.  We assume that the 
set $M$ is a minimal system of generators for the semigroup $\Gamma(M)$. 
The semigroup is known as a semigroup generated by concatanation of two 
arithmetic sequences. Given $p\geq 3$, what conditions on $a,b,d$ ensure 
that $\beta_{i}(\mathfrak{p}(M))$ are unbounded for every $i$?
\end{enumerate}

\section{The Derivation Module}
Let $k$ be a field of characteristic zero. We now consider the algebroid case 
and denote the coordinate ring of the algebroid monomial curve defined by the 
integers $n_0,\ldots, n_p$; $p\geq 2$, by the same notation $A$;
$$A=k[[x_0,\,\ldots,\, x_p]]/\frak{p}(n_0,\ldots, n_p)=k[[t^{n_{0}}, \ldots ,t^{n_{p}}]].$$
The module of $k$-derivations 
${\rm Der}_{k}(A)$, known as the \textit{Derivation module of $A$} is 
an extremely important object of study. A good understanding of the 
derivation module is essential for various reasons and surely from the 
perspective of the Zariski-Lipman conjecture, which says that over a 
field of characteristic $0$, $A$ is a polynomial ring over $k$ if 
${\rm Der}_{k}(A)$ is a free $A$-module; see \cite{hochster}, 
\cite{mullerpatil}. An explicit generating set for the derivation module 
has been constructed by Patil and Singh \cite{patsingh} under the assumption 
that the integers $n_0,\ldots, n_p$ form a minimal almost arithmetic sequence. 
There are discussions on other curves in \cite{patil} which has interesting 
results on their derivation modules. In \cite{patseng}, the minimal number 
of generators for the derivation modules has been calculated. This also gives 
a formula for the last Betti number for the defining ideal 
$\frak{p}(n_0,\ldots, n_p)$. The following theorem was proved by Kraft \cite{kraft}; page 
875, which reduces the computation of $\mu ({\rm Der}_{k}(A))$ to a counting problem; 
see \cite{patseng}.
\medskip 

\begin{theorem} 
The set $\left\{ t^{\alpha + 1} \frac{d}{dt} \mid \alpha \in
\Delta' \cup \{0\}\right\}$ is a minimal set of generators 
for the $A$-module ${\rm Der}_{k}(A)$, where $\Gamma(n_0,\ldots, n_p)$ is 
the numerical semigroup generated by the sequence $n_0,\ldots, n_p$ 
of positive integers, 
$$\Gamma(n_0,\ldots, n_p)_{+} := 
\Gamma(n_0,\ldots, n_p) \setminus \{0\};$$
$$\Delta := \{ \alpha \in \mathbb{Z}^{+} \mid \alpha + \Gamma(n_0,\ldots, n_p)_{+} \subseteq    
\Gamma(n_0,\ldots, n_p)\};$$
$$\Delta' := \Delta\setminus \Gamma(n_0,\ldots, n_p).$$
In particular, $\mu ({\rm Der}_{k}(A)) = {\rm card}\left(\Delta'\right) +1$. 
\end{theorem}
\medskip

\begin{enumerate}
\item[Q6.] Compute $\mu ({\rm Der}_{k}(A))$, where $k$ is a field of 
characteristic $0$ and $A = k[[x_0,\,\ldots,\, x_p]]/\frak{p}(n_0,\ldots, n_p)$.
\end{enumerate}

\section{Smoothness of blowups and Francia's conjecture}
Let $R$ be a regular local ring. Let $\mathfrak{p}$ be an ideal 
in $R$. The Rees algebra $\mathcal{R}(\mathfrak{p}) = R[\mathfrak{p}T]$ 
is called the Rees Algebra of the ideal $\mathfrak{p}$. The Blowup 
Algebras, in particular the Rees algebra and the associated graded ring 
$\mathcal{G}(\mathfrak{p}) = 
\mathcal{R}(\mathfrak{p})/\mathfrak{p}\mathcal{R}(\mathfrak{p})$ 
of the ideal $\mathfrak{p}$ 
play a crucial role in the birational study of curves. See \cite{vasco} 
for learning the basics and related results on blowup algebras. The projective 
scheme ${\rm Proj}(\mathcal{R}(\mathfrak{p}))$ is the 
blowup of ${\rm Spec}(R)$ along $V(\mathfrak{p})$. Francia's conjecture 
is the following:
\medskip

\noindent\textbf{Francia's Conjecture.} If $\mathfrak{p}$ is a 
prime ideal in a regular local ring $R$ with ${\rm dim}(R/\mathfrak{p})=1$ and 
if ${\rm Proj}(\mathcal{R}(\mathfrak{p}))$ is a smooth projective scheme then 
$\mathfrak{p}$ is a complete intersection. 
\medskip

\noindent The conjecture stated in the context of monomial curves is the following:
\medskip

\noindent\textbf{Francia's Conjecture for monomial curves.} If ${\rm Proj}(\mathcal{R}(\mathfrak{p}(n_{0}, \ldots, n_{p})))$ is a smooth projective scheme, then $\mu(\mathfrak{p}(n_{0}, \ldots, n_{p}))=p$, 
in other words $\mathfrak{p}(n_{0}, \ldots, n_{p})$ is a complete intersection.
\medskip

The conjecture is false in general. It was proved in \cite{jhonmorey}, 
that, if $R=\mathbb{Q}[x,y,z]$ and if $\mathfrak{p}$ denotes the ideal 
generated by the maximal minors of the matrix
$$
\left(\begin{matrix}
0 & x-y & y-z\\
z-y & 0 & x-y+z\\
z & 0 & x-2y\\
x & y & z
\end{matrix}\right),
$$
then $\mathfrak{p}$ is a prime ideal of codimension $2$ and $R[\mathfrak{p}t]$ 
is smooth. This produces a counter example for Francia's conjecture over the 
base filed $k=\mathbb{Q}$. However, over an algebraically closed field this 
ideal does not remain a prime ideal. Therefore, the conjecture is open in 
its geometric form, that is, over $k=\mathbb{C}$. The conjecture has been proved 
for certain cases in \cite{carrollvalla}. It seems that the hardest 
part in solving this conjecture is computing the equations defining the Rees 
algebra as an affine algebra over $R$; see \cite{mukhsen2}. In \cite{mukhsen1} 
we have proved that the conjecture is true if the integers $n_{0}, \ldots, n_{p}$ 
form an arithmetic sequence. We could do this without explicitly computing all 
the equations defining the Rees algebra, in order to show that the projective 
scheme ${\rm Proj}(\mathcal{R}(\mathfrak{p}(n_{0}, \ldots,n_{p})))$ is not 
smooth. This answered the conjecture in affirmative because we knew that 
$\mathfrak{p}(n_{0}, \ldots, n_{p})$ could never be a complete intersection 
if $p\geq 3$ and $n_{0}, \ldots, n_{p}$ form an arithmetic sequence; see 
\cite{malooseng}. One case which has not been studied so far is when 
the integers $n_{0}, \ldots, n_{p}$ form an almost arithmetic sequence, 
that is, all but one of these integers form an arithmetic sequence.
\medskip

\begin{enumerate}
\item[Q 7.] Answer Francia's conjecture when $n_{0}, \ldots, n_{p}$ 
form a minimal almost arithmetic sequence of integers.
\end{enumerate}

\section{The Set theoretic complete intersection conjecture}
To quote M. Hochster \cite{hochsterproblems}, ``This seems incredibly difficult". Let us 
state the problem in the context of affine monomial curves. 
\begin{enumerate}
\item[Q6.] There exist $p$ polynomials $f_{1}, \ldots , f_{p}$ 
in $k[x_{0}, \ldots , x_{p}]$, such that 
$$\frak{p}(n_0,\ldots, n_p)=\sqrt{(f_{1}, \ldots , f_{p})}.$$
\end{enumerate}
The conjecture is true if charactertistic of $k$ is a prime, by the 
famous result proved by Cowsik and Nori \cite{cowsiknori}. However, 
very little is known in characteric zero, except for some special cases. 
We refer to \cite{lyubeznik} for a detailed account of research done on this 
conjecture. Patil \cite{patil1} proved that affine monomial curves defined by an almost 
arithmetic sequence is a set theoretic complete intersection. Bresinsky 
proved the conjecture in dimension 3 in \cite{bre4} and that in dimension 
$4$ with symmetry condition in \cite{bre5}. 
These results were subsequently generalized by Thoma \cite{thoma}. Eto  \cite{eto} proved the 
conjecture when the integers define a balanced semigroup in dimension 
$4$. The general question is still open. 
\medskip

An important technique using Gr\"{o}bner basis was mentioned in \cite{robvalla}, 
which says the following: If $V$ is a projective variety of codimension $r$ in $\mathbb{P}^{n}$, 
which is the projective closure of an affine variety $W$ in $\mathbb{A}^{n}$, and if 
$I(W) =  \sqrt{(f_{1}, \ldots , f_{r})}$, such that $\{f_{1}, \ldots , f_{r}\}$ is a 
Gr\"obner basis for the ideal generated by $\{f_{1}, \ldots , f_{r}\}$, then 
$I(V) = I(W)^{h} = \sqrt{(f_{1}^{h}, \ldots , f_{r}^{h})}$. The authors have proved 
that the rational normal curve in $\mathbb{P}^{n}$ is a set-theoretic complete 
intersection by the technique mentioned above. The paper \cite{valla} precedes 
\cite{robvalla} and proves 
several interesting results on ideals defined by determinantal conditions. The 
set-theoretic complete intersection property of the rational normal curve forms 
an important step towards proving the set-theoretic complete intersection property 
of an affine monomial curve defined by an almost arithmetic sequence in \cite{patil1}. 
\medskip

An interesting observation is due to Cowsik \cite{cowsik}, which is the following: 
For an one dimensional prime $\mathfrak{p}$ in a regular local ring $R$, if the 
symbolic Rees algebra $\oplus_{n=1}^{\infty}\mathfrak{p}^{(n)}$ is Noetherian then 
$\mathfrak{p}$ is a set-theoretic complete intersection. However, it turns 
out that the Noetherian property of the symbolic Rees algebra or equivalently the 
finite generation of the symbolic Rees algebra does not happen commonly. 
Huneke \cite{hun1} proved that $\oplus_{n=1}^{\infty}\mathfrak{p}(3, n_{1}, n_{2})^{(n)}$ 
is Noetherian if $3< n_{1} < n_{2}$. In a latter paper \cite{hun2}, Huneke gave 
a useful necessary and sufficient condition for the Noetherian property 
of the symbolic Rees algebra for a one dimensional prime in a regular local ring of 
dimension $3$ with infinite residue field. Some more results pertaining to the finite 
generation of the symbolic Rees algebra of monomial primes can be found in 
\cite{cutkosky} and \cite{hema}.

\bibliographystyle{amsalpha}

\begin{thebibliography}{A}

\bibitem{alfonsin}
J.L. Ramirez Alfonsin, \emph{The Diophantine Frobenius problem}, Oxford University Press, 
NY, 2005.

\bibitem{bh}
W. Bruns, J. Herzog, \emph{Cohen Macaulay Rings}, Cambridge University Press, NY, 1993.

\bibitem{bre}
 H. Bresinsky, \emph{On Prime Ideals with Generic Zero $x_i=t^{n_{i}}$}, Proceedings of the American Mathematical Society, 47(2)(1975).

\bibitem{bre1} H. Bresinsky, \emph{Symmetric semigroups of integers generated by $4$ elements}, 
Manuscripta Math. 17(1975), 205--219.

\bibitem{bre2} H. Bresinsky, \emph{Monomial Gorenstein Ideals}, 
Manuscripta Math. 29(1979), 159--181.

\bibitem{bre3} H. Bresinsky, \emph{Binomial generating sets for monomial curves 
with applications in $\mathbb{A}^{4}$}, Rend. Sem. Mat. Univers. Politecn. Torino, 
46(3)(1988).

\bibitem{bre4} H. Bresinsky, \emph{Monomial space curves in $\mathbb{A}^{3}$ as set-theoretic complete intersections.} Proc. Amer. Math. Soc. 75 (1979), no. 1, 23-24. 

\bibitem{bre5} H. Bresinsky, \emph{Monomial Gorenstein curves in 
$\mathbb{A}^{4}$ as set-theoretic complete intersections.} Manuscripta Math. 
27(4)(1979), 353--358.

\bibitem{brehoa} 
H. Bresinsky, L. T. Hoa, \emph{Minimal generating sets for a family of monomial curves 
in $\mathbb{A}^{4}$}, Commutative algebra and algebraic geometry (Ferrara), Lecture 
Notes in Pure and Appl. Math. 206 (Dekker, New York, 1999) 5--14.

\bibitem{carrollvalla} L. O’Carroll, G. Valla, \emph{On the smoothness of blowups}, 
Commun. Algebra 25(6)(1997) 1861-–1872.

\bibitem{cowsik} R. C. Cowsik, \emph{Symbolic powers and the number of defining equations.} 
Algebra and its applications (New Delhi, 1981), 13-14, Lecture Notes in Pure and Appl. Math., 
91, New York, 1984.

\bibitem{cowsiknori} R. C. Cowsik, M. V. Nori, \emph{Affine curves 
in characteristic $p$ are set theoretic complete intersections.} 
Invent. Math. 45(1978), no. 2, 111-114.

\bibitem{clo}
D. Cox, J. Little, and D. O'Shea, \emph{Ideals, varieties, and algorithms: an introduction to computational algebraic geometry and commutative algebra}, New York: Springer,(2007).

\bibitem{cox}
D. Cox, J. Little and D. O'Shea, \emph{Using Algebraic Geometry}, 
GTM 185 (Springer-Verlag New York, 1998).

\bibitem{cutkosky} S.D. Cutkosky, \emph{Symbolic algebras of monomial primes}, 
J. Reine Angew. Math. 416: 71--89 (1991).

\bibitem{dl} Decker, W., Lossen, C., \emph{Computing in
Algebraic Geometry A Quick Start using SINGULAR}, Hindustan Book Agency, 
October 2005.

\bibitem{eto} K. Eto, \emph{The monomial curves associated with balanced semigroups are set-theoretic complete intersections}, Journal of Algebra, 
319 (2008), 1355--1367.

\bibitem{GAP4} The GAP~Group, {\em GAP -- Groups, Algorithms, and Programming, 
Version 4.8.6}; 2016. \newblock {https://www.gap-system.org} (2016). 
 
\bibitem{gss1}
P. Gimenez, I. Sengupta and H. Srinivasan, 
\emph{Minimal free resolution for certain affine monomial curves}. 
In: {\em Commutative Algebra and its Connections to Geometry} 
({\em PASI 2009}), A.~Corso and C.~Polini Eds, Contemp. Math. 555 
(Amer.Math. Soc., 2011) 87--95.

\bibitem{gss2} P. Gimenez, I. Sengupta, H. Srinivasan, \emph{Minimal graded free resolutions for monomial curves defined by arithmetic sequences}, J.Algebra  388(2013), 294--310.

\bibitem{singular} 
Decker, W.; Greuel, G.-M.; Pfister, G.; Sch{\"o}nemann, H., 
\emph{\newblock {\sc Singular} {4.0.1} --- {A} computer algebra system for polynomial computations}. \newblock {http://www.singular.uni-kl.de} (2014).

\bibitem{singularbook} 
Greuel, G.-M.; Pfister, G. \emph{A Singular introduction to Commutative Algebra}, 
Springer-Verlag, 2002.

\bibitem{herzog}
J. Herzog, \emph{Generators and relations of abelian semigroups and semigroup rings}, 
Manuscripta Mathematica, Vol. 2, No. 3, (1970), 175--193.

\bibitem{herstamate} J. Herzog, Dumitru I. Stamate, \emph{On the defining equations 
of the tangent cone of a numerical semigroup ring}, Journal of Algebra, 418(2014), 8--28.

\bibitem{hochster} M. Hochster, \emph{The Zariski-Lipman conjecture in 
the graded case}, Journal of Algebra 47(1977), 411--424.

\bibitem{hochsterproblems} M. Hochster, \emph{Thirteen Open Questions in Commutative Algebra}. 
Available at http://www.math.lsa.umich.edu/~hochster/Lip.text.pdf

\bibitem{hun1} C. Huneke, \emph{On the finite generation of symbolic blowups}, 
Math.Z 179(1982), 465--572.

\bibitem{hun2} C. Huneke, \emph{Hilbert Functions and Symbolic Powers}, 
Michigan Math.J. 34(1987), 293--318.

\bibitem{iyengar}
S. B. Iyengar, G. J. Leuschke, A. Leykin,  C. Miller, E. Miller, A. K. Singh, U. Walther, \emph{Twenty-Four Hours of Local Cohomology}, American Mathematical Society, 2007.

\bibitem{jhonmorey} M.R. Johnson, S. Morey, \emph{Normal blow-ups and their expected 
defining equations}, J. Pure Appl. Algebra 215(2001) 162, 303-–313.

\bibitem{kraft} J. Kraft, \emph{Singularity of monomial curves in $\mathbb{A}^{3}$ and
Gorenstein monomial curves in $\mathbb{A}^{4}$}, Canadian J. Math. 37(1985) 872--892.

\bibitem{lyubeznik} G. Lyubeznik, \emph{A survey of problems and results on the
number of defining equations}, Commutative algebra
(Berkeley, CA, 1987), 375-390, MSRI Publ., 15, Springer, NY, 1989. 

\bibitem{malooseng}
A. K. Maloo and I. Sengupta, \emph{Criterion for complete intersection of certain
monomial curves}. In: {\em Advances in Algebra and Geometry} ({\em Hyderabad, 2001}) (Hindustan Book Agency, New Delhi, 2003) 179--184.

\bibitem{matsumura}
H. Matsumura, \emph{Commutative Ring Theory}, Cambridge University Press, NY, 1986.

\bibitem{mss1} R. Mehta, J. Saha, I. Sengupta, \emph{Betti numbers of 
Bresinsky's curves in  $\mathbb{A}^4$}. Journal of Algebra and Its Applications
Vol. 18, No. 8 (2019) 1950143 (14 pages). 

\bibitem{mss2} R. Mehta, J. Saha, I. Sengupta, \emph{Numerical Semigroups generated by concatenation of arithmetic sequences}; to appear in 
Journal of Algebra and Its Applications. DOI: 10.1142/S0219498821501620

\bibitem{mukhsen1} D. Mukhopadhyay, I. Sengupta, \emph{On the Smoothness of Blowups 
for Certain Monomial Curves}, Beitr\"{a}ge zur Algebra und Geometrie, 53(2012), 89--95, 
Springer. 

\bibitem{mukhsen2} D. Mukhopadhyay, I. Sengupta, \emph{The Rees Algebra for Certain 
Monomial Curves}, Ramanujan Mathematical Society-Lecture Notes Series 17 (2003), Proc. 
CAAG 2010, pp.199--218.

\bibitem{mullerpatil} G. M\"{u}ller, D.P. Patil, \emph{The 
Herzog-Vasconcelos conjecture for affine semigroup rings} 
Communications in Algebra, 27(7)(1999) 3197--3200.

\bibitem{patil} 
D.P. Patil, \emph{Generators for the derivation modules and the defining ideals of certain affine curves}, \textit{Thesis}, TIFR-Bombay University, 1989.

\bibitem{patil1}
{D. P. Patil, Certain monomial curves are set theoretic complete intersections,
{\em Manuscripta Math.} 68 (1990) 399--404.
}

\bibitem{patsingh} D.P. Patil, B. Singh, \emph{Generators for the
derivation modules and the relation ideals of certain curves}, 
Manuscripta Math. 68(1990) 327-335.

\bibitem{patil2}
D.P. Patil, \emph{Minimal sets of generators for the relation ideals 
of certain monomial curves}, Manuscripta Math. 80 (1993) 239--248.

\bibitem{peeva}
I. Peeva, \emph{Graded Syzygies}, Springer-Verlag London Limited, 2011.

\bibitem{robvalla} L. Robbiano, G. Valla, \emph{On set-theoretic complete intersections 
in the projective space}, Rendiconti del Seminario Matematico e Fisico di Milano, 
Vol. LIII (1983).
 
\bibitem{rgs} 
J.C. Rosales, P.A. Garc\'{i}a-S\'{a}nchez, \emph{Numerical Semigroups}, 
Springer, (2009).

\bibitem{ros} J.C.Rosales, \emph{Symmetric Numerical Semigroups with Arbitrary Multiplicity and Embedding Dimension}, Proceedings of the American Mathematical Society, 29(8)(2001), 2197--2203.

\bibitem{rst} A.K. Roy, I. Sengupta, G. Tripathi, \emph{Minimal graded free resolutions for monomial curves defined by almost arithmetic sequences}, Communications in Algebra 
45(2), 2017, pp. 521-551.

\bibitem{patseng} D.P. Patil, I. Sengupta, \emph{Minimal set of generators for 
the derivation module of certain affine monomial curves}, Comm. Algebra 27(11)(1999) 
5619-5631.

\bibitem{seng1} I. Sengupta, \emph{A Gr\"{o}bner basis for certain affine monomial curves}, 
Comm. Algebra 31(3) (2003) 1113--1129.

\bibitem{seng2} I. Sengupta, \emph{A minimal free resolution for certain monomial curves 
in $\mathbb{A}^{4}$}, Comm. Algebra 31(6) (2003) 2791--2809.

\bibitem{stamate} Dumitru I. Stamate, \emph{Betti numbers for numerical 
semigroup rings}, arxiv:1801.00153v1[math.AC], 2017.

\bibitem{hema} Hema Srinivasan, \emph{On Finite Generation of Symbolic Generation of 
Monomial Primes}, Comm. Alg. 19(9) 2557--2564 (1991).

\bibitem{ian} I. Stewart, \emph{How hard can it be?}, New Scientist, 21 June 2008, 46--49.

\bibitem{sturmfels} B. Sturmfels, \emph{Gr\"{o}bner Bases and Convex Polytopes}, University Lecture Series, Volume 8, AMS, 1996.

\bibitem{thoma} A. Thoma, \emph{On the set-theoretic complete intersection problem 
for monomial curves in $\mathbb{A}^{n}$ and $\mathbb{P}^{n}$}, Journal of Pure and 
Applied Algebra 104 (1995) 333-344.

\bibitem{valla} G. Valla, \emph{On determinantal ideals which are set-theoretic complete intersections}, 
Compositio Mathematica 42(1) 3-11 (1981).

\bibitem{vasco} W. Vasconcelos, \emph{Arithmetic of Blowup Algebras}, 
LMS Lecture Note Series, vol. 195. Cambridge University Press, UK, 1994.
\end{thebibliography}

\end{document}